\documentclass[12 pt]{article}
\usepackage{array,color}
\usepackage{amsbsy}
\usepackage{amsmath}
\usepackage{amssymb}
\usepackage{amsfonts}
\usepackage{latexsym}
\usepackage{epsf}
\usepackage{epsfig}
\usepackage{graphicx}
\usepackage{fancybox,fancyhdr,times,subfigure}
\usepackage{epsfig, color, graphicx, rotate, colordvi}
\usepackage{a4wide}
\usepackage[numbers]{natbib}
\usepackage{mathrsfs,dsfont,amssymb}
\usepackage{bm}

\newcommand{\bbI}{{\bf I}}

\newcommand{\bbF}{{\bf F}}
\newcommand{\bbe}{{\bf e}}

\newcommand{\bbR}{{\bf R}}
\newcommand{\bbgL}{{\bf \Lambda}}
\newcommand{\bgma}{{\bf \gamma}}
\newcommand{\bbM}{{\bf M}}
\newcommand{\bPhi}{{\bf \Phi}}

\newcommand{\bqn}{\begin{eqnarray*}}
\newcommand{\eqn}{\end{eqnarray*}}
\newcommand{\rE}{{\textrm{E}}}

\newcommand{\bqa}{\begin{eqnarray}}
\newcommand{\eqa}{\end{eqnarray}}

\newcommand{\non}{\nonumber\\}
\newcommand{\ep}{\varepsilon}
\newcommand{\rP}{\textrm{P}}

\newcommand{\la}{\lambda}
\newcommand{\proof}{\textsc{Proof}}

\newtheorem{thm}{Theorem}[section]

\newtheorem{remark}{{\it Remark}}[section]
\newtheorem{lemma}{Lemma}[section]
\begin{document}

\begin{center}
A Note on the Limiting Spectral Distribution of a Symmetrized Auto-Cross Covariance Matrix\\[2ex]
Zhidong Bai\footnote{The research of this author was supported by NSF China 11171057, the Program for Changjiang Scholars and the Innovative Research Team in University, the Fundamental Research Funds for the Central Universities and the NUS Grant R-155-000-141-112.}  and Chen Wang
\end{center}
\begin{abstract}
In Jin et al. (2014), the limiting spectral distribution (LSD) of a symmetrized auto-cross covariance matrix is derived using matrix manipulation, with finite $(2+\delta)$-th moment assumption. Here we give an alternative method using a result in Bai and Silverstein (2010), in which a weaker condition of finite 2nd moment is assumed.
\end{abstract}
\section{Introduction}
Consider a large dimensional dynamic $k$-factor model with lag $q$ taking the form of  $$\bbR_{t}=\sum_{i=0}^{q}\bbgL_{i}
\bbF_{t-i}+\bbe_{t}, \quad t=1,...,T$$
where \textbf{$\bbgL_{i}$}'s are $N\times k$ non-random matrices with full rank. For $t=1,...,T$,
$\bbF_{t}$'s are $k$-dimensional vectors of independent identically distributed (i.i.d.) standard complex components and $\bbe_{t}$'s are $N$-dimensional vectors of i.i.d. complex components with mean zero and finite second moment $\sigma^2$, independent of $\bbF_{t}$. This can also be considered as
a type of \emph{information-plus-noise model} (Dozier and Silverstein, 2007a, b; Bai and Silverstein, 2012) where the information comes from the summation part and the noise is $\bbe_t$'s. Here both $k$ and $q$ are fixed but unknown, while both $N$ and $T$ tend to $\infty$ proportionally.\\

Under this high dimensional setting, an important statistical problem is the estimation of $k$ and $q$ (Bai and Ng, 2002; Harding, 2012).
To this objective, the following two variables are defined for fixed non-negative integer $\tau$, namely:
$$\bPhi_N(\tau)=\frac{1}{2T}\sum_{j=1}^{T}(\bbR_{j}\bbR_{j+\tau}^{*}+\bbR_{j+\tau}\bbR_{j}^{*})$$ and
$$\bbM_N(\tau)=\sum_{j=1}^{T}(\bgma_{j}\bgma_{j+\tau}^{*}+\bgma_{j+\tau}\bgma_{j}^{*}),$$ where $\bgma_{j}=\frac{1}{\sqrt {2T}}\bbe_{j}$ and $*$ denotes the conjugate transpose.\\

Note that when $\tau=0$, we have $\bbM_N(\tau)=\frac1T\sum_{j=1}^{T}\bbe_{j}\bbe_{j}^{*}$, which is a sample covariance matrix, whose LSD follows MP law (Mar\v{c}enko and Pastur, 1967) with density $$f_c(x)=\frac1{2\pi cx}\sqrt{(b_c-x)(x-a_c)}, x\in[a_c,b_c]$$
and a point mass $1-1/c$ at the origin if $c>1$. Here $c=\lim_{N\to\infty} N/T$, $a_c=(1-\sqrt{c})^{2}$ and $b_c=(1+\sqrt{c})^{2}$.\\

Moreover, if we write $$\boldsymbol{\bbgL}=(\boldsymbol{\bbgL}_0,\boldsymbol{\bbgL}_1,\cdots,\boldsymbol{\bbgL}_q)_{N\times k(q+1)},$$
then the covariance matrix of $\bbR_{t}$ will be similar to
$$\begin{pmatrix}
\sigma^2 \bbI+\boldsymbol{\bbgL}^{*}\boldsymbol{\bbgL} & \textbf{0}\\
\textbf{0} & \sigma^2 \bbI
\end{pmatrix},$$
with the size of the upper block and lower block $k(q+1)$ and $N-k(q+1)$, respectively. Thus, we have a \emph{spiked population model} (Johnstone, 2001; Baik and Silverstein, 2006; Bai and Yao, 2008).
In fact, under certain conditions, $k(q+1)$ can be estimated by counting the number of eigenvalues of $\bPhi_N(0)$ that go beyond the certain phase trasition point. Therefore, it remains to estimate one of $k$ and $q$. To this end, it is necessary to investigate the LSD of $\bbM_N(\tau)$ for at least one $\tau \ge 1$. As such, Jin et al. (2014) has established the following result.
\begin{thm} (Theorem 1.1 in Jin et al. (2014))\label{thm1} Assume:
 \begin{description}
   \item(a) $\tau \ge1$ is a fixed integer.
   \item(b) $\bbe_{k}=(\ep_{1k},\cdots,\ep_{Nk})'$, $k=1,2,...,T+\tau$, are $N $ dimensional vectors of independent standard complex components with $\sup_{1\le i\le N,1\le t\le T+\tau}\rE|\ep_{it}|^{2+\delta}\le M <\infty$ for some $\delta\in(0,2]$, and for any $\eta>0$, \bqa \label{trun1} \frac{1}{\eta^{2+\delta}N T}\sum_{i=1}^N\sum_{t=1}^{T+\tau}\rE(|\ep_{it}|^{2+\delta}I(|\ep_{it}|\geq \eta T^{1/(2+\delta)}))=o(1).
\eqa
   \item(c) $N/(T+\tau)\rightarrow c>0$ as $N,T\to \infty$.
   \item(d) $\bbM_N=\sum_{k=1}^{T}(\bgma_{k}\bgma_{k+\tau}^{*}+\bgma_{k+\tau}\bgma_{k}^{*}),$ where $\bgma_{k}=\frac{1}{\sqrt {2T}}\bbe_{k}$.
 \end{description}
 Then as $N,T\rightarrow\infty$, $F^{\bbM_N}\overset {D}{\rightarrow} F_c$ a.s. and $F_c$ has a density function given by
 \begin{center}
 $\phi_c(x)=\frac{1}{2c\pi}\sqrt{\frac{y_0^2}{1+y_0}-(\frac{1-c}{|x|}+\frac{1}{\sqrt{1+y_0}})^2}$, $ |x|\leq a,$
 \end{center}
 where
\begin{eqnarray*}
a=\left\{\begin{array}{cc}\frac{(1-c)\sqrt{1+y_1}}{y_1-1},& c \neq 1,\\
2,& c=1,
\end{array}\right.
\end{eqnarray*}
$y_0$ is the largest real root of the equation: $y^3-\frac{(1-c)^2-x^2}{x^2}y^2-\frac{4}{x^2}y-\frac{4}{x^2}=0$ and $y_1$ is the only real root of the equation: \bqa((1-c)^2-1)y^3+y^2+y-1=0\eqa such that $y_1>1$ if $c<1$ and $y_1\in (0,1)$ if $c>1$.
Further, if $c>1$, then $F_c$ has a point mass $1-1/c$ at the origin.
\end{thm}
In Jin et al. (2014), the key step of the proof of Theorem \ref{thm1} is to establish that the Stieltjes transform $m$ of $F_c$ satisfies \bqa\label{eq1}
(1-c^2m^2(z))(c+czm(z)-1)^2=1,
\eqa
from which four roots are obtained:
\bqn
m_1(z)&=&\frac{(\frac{1-c}{z}+\sqrt{1+y_0})+\sqrt{(\frac{1-c}{z}-\frac{1}{\sqrt{1+y_0}})^2-\frac{y_0^2}{1+y_0}}}{2c}\\
m_2(z)&=&\frac{(\frac{1-c}{z}+\sqrt{1+y_0})-\sqrt{(\frac{1-c}{z}-\frac{1}{\sqrt{1+y_0}})^2-\frac{y_0^2}{1+y_0}}}{2c}\\
m_3(z)&=&\frac{(\frac{1-c}{z}-\sqrt{1+y_0})+\sqrt{(\frac{1-c}{z}+\frac{1}{\sqrt{1+y_0}})^2-\frac{y_0^2}{1+y_0}}}{2c}\\
m_4(z)&=&\frac{(\frac{1-c}{z}-\sqrt{1+y_0})-\sqrt{(\frac{1-c}{z}+\frac{1}{\sqrt{1+y_0}})^2-\frac{y_0^2}{1+y_0}}}{2c}.\\
\eqn
Here $y_0$ is the largest real root of the equation:
$$f(y):=y^3-\frac{(1-c)^2-z^2}{z^2}y^2-\frac{4}{z^2}y-\frac{4}{z^2}=0.$$
Note that all the three roots of $f(y)=0$ give the same set of $m_i$'s, up to a permutation order,
and our choice $y_0$ as the largest real root is only for the sake of simplicity.\\

For the four $m_i$'s, after some justification, we have
\begin{eqnarray*}
m(z)=\left\{\begin{array}{cc}m_1(z),& z< 0,\\
m_3(z),& z>0.
\end{array}\right.
\end{eqnarray*}
The density function is then derived using the inversion formula of the Stieltjes transform.\\

Figures 1 and 2 display the density functions $\phi_c(x)$ with $c<1$ and $c>1$, respectively.
From these two figures, it is shown that as $c$ increases, the support of $\phi_c(x)$ gets wider, and $\phi_c(x)$ achieves the maximum at $x=0$ which is sharper as $c$ gets closer to 1.\\

\begin{figure}[!ht]
\begin{center}
 \includegraphics[scale=0.4]{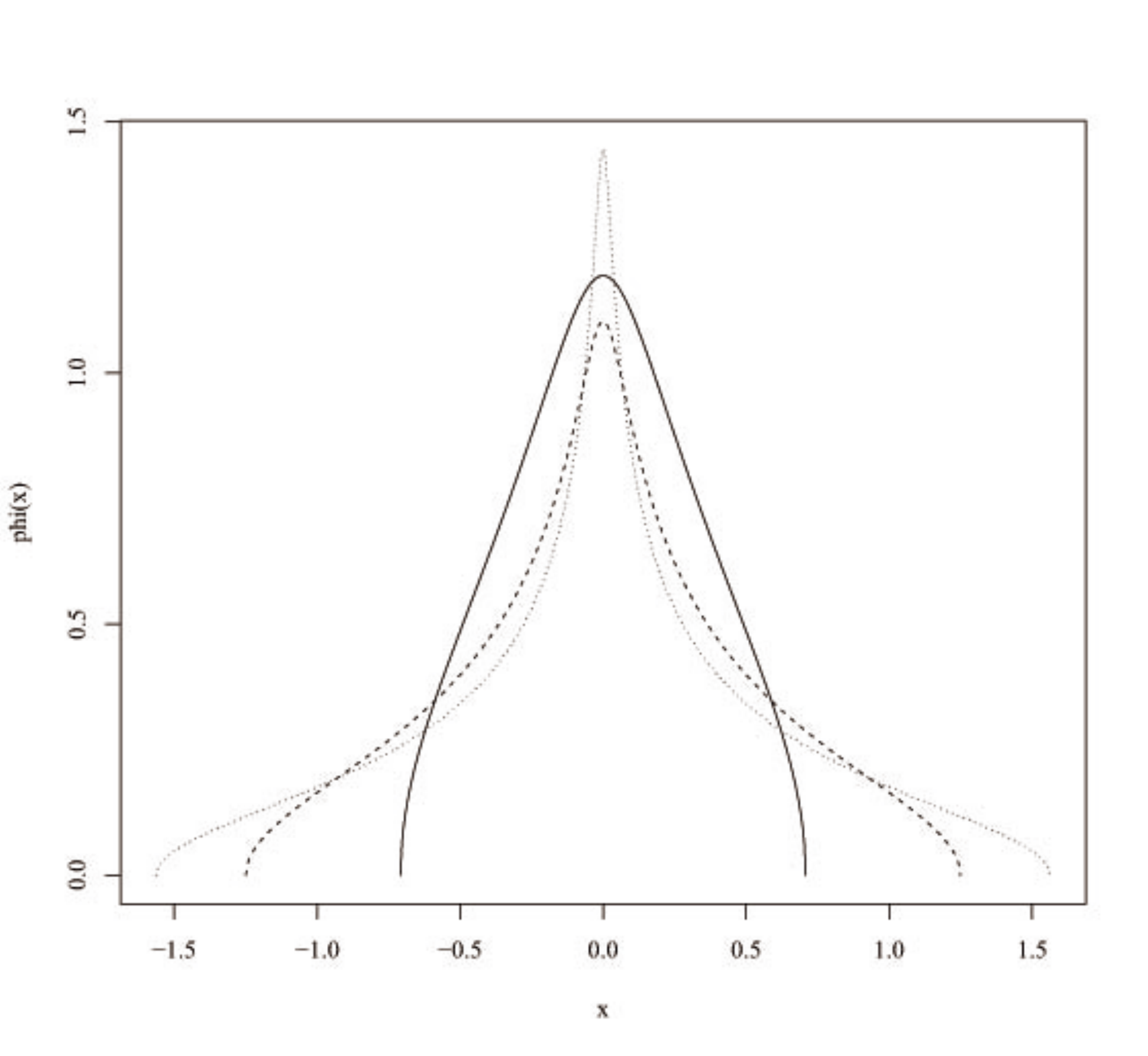}
 \caption{Density functions $\phi_c(x)$ of the LSD of $\bbM_N$ with $c=0.2$ (the solid  line), $c=0.5$ (the dashed line) and $c=0.7$ (the dotted line).}\label{fig1}
\end{center}
\end{figure}

\begin{figure}[!ht]
\begin{center}
 \includegraphics[scale=0.4]{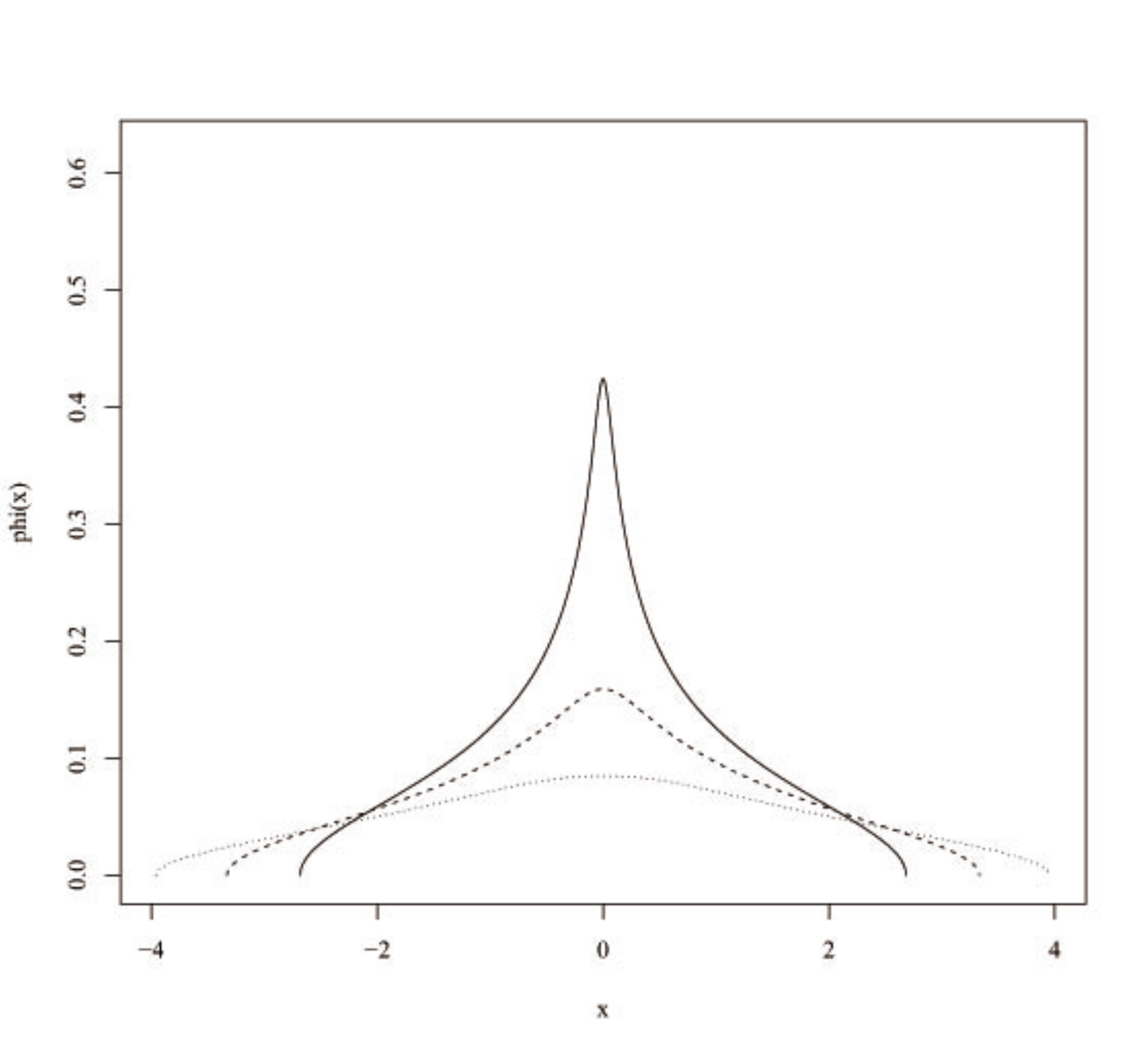}
 \caption{Density functions $\phi_c(x)$ of the LSD of $\bbM_N$ with $c=1.5$ (the solid line), $c=2$ (the dashed line) and $c=2.5$ (the dotted line). Note that the area under each density function curve is $1/c$.}\label{fig2}
\end{center}
\end{figure}
The goal of the paper is devoted to giving a more direct method of deriving (\ref{eq1}), by using Theorem \ref{thm2} in Bai and Silverstein (2010). It is worth noting that for our method to work, we only require the finiteness of the 2nd moment of the underlying random variable, which is weaker than the finite $(2+\delta)$-th moment requirement in Jin et al. (2014). Once (\ref{eq1}) is obtained, Theorem \ref{thm1} will follow by employing the technique in Jin et al. (2014) and thus will not be presented again.
\section{Notation}
Before proceeding, it is necessary to rewrite $\bbM_N(\tau)$ into another form. For any $\tau\ge1$ fixed, write
\bqn
&&\bbM_N(\tau)\\
&=&\sum_{k=1}^{T}(\bgma_{k}\bgma_{k+\tau}^{*}+\bgma_{k+\tau}\bgma_{k}^{*})\\
&=&\frac1{2T}\sum_{k=1}^{T}(\bbe_{k}\bbe_{k+\tau}^{*}+\bbe_{k+\tau}\bbe_{k}^{*})\\
&=&\frac1{T}(\bbe_1,\bbe_2,\cdots,\bbe_{T+\tau-1},\bbe_{T+\tau})\begin{pmatrix}
  0 & \cdots &  \frac12 & \cdots & 0 \\
  \vdots & \ddots & 0 &  \frac12 &  \vdots\\
 \frac12 &  0 &\ddots&0 & \frac12  \\
  \vdots & \frac12& 0 & \ddots & \vdots   \\
  0 &\cdots & \frac12& \cdots & 0  \\
 \end{pmatrix}\begin{pmatrix}
\bbe_1^* \\
 \bbe_2^*\\
 \vdots\\
 \bbe_{T+\tau-1}^*\\
  \bbe_{T+\tau}^*
 \end{pmatrix}\\
 &\equiv&\frac1{T}\textbf{X}_T\textbf{C}_{T,\tau}\textbf{X}_T^*,
 \eqn
where the two bands of $\frac12$'s are $\tau-$distance from the main diagonal.
\section{A Useful Lemma}
\begin{lemma}\label{key} As $n\to\infty$, the empirical spectral distribution (ESD) of $\textbf{C}_{n,\tau}$ tends to $H$, which is an Arcsine distribution with density function
\begin{center}
 $H'(t)=\frac{1}{\pi\sqrt{1-t^2}}, \quad t\in(-1,1).$
 \end{center}
\end{lemma}
\proof. Fix $n$ and let $\la$ be an eigenvalue of $\textbf{C}_{n,\tau}$.\\
Define $D_n=D_{n,\tau}=\det(\la\textbf{I}-\textbf{C}_{n,\tau})$ (for simplicity, we omit $\tau$ from the subscript).\\
When $n<\tau$, all the entries of $\textbf{C}_{n,\tau}$ are 0 and hence we have $D_n=\la^n$.
For $n\ge\tau$, expand along the first row, and we have $D_n=\la D_{n-1}+\frac{(-1)^\tau}2\tilde D_{n-1}.$\\
Expand along the first column of the matrix wrt $\tilde D_{n-1}$, and we have $\tilde D_{n-1}=\frac{(-1)^{\tau-1}}2D_{n-2}.$\\
Therefore, for $n\ge\tau$, we have
\bqn
D_n=\la D_{n-1}-\frac14D_{n-2}
\eqn
Solve the characteristic equation $x^2=\la x-\frac14$ and we have $\la_{1,2}=\frac{\la\pm\sqrt{\la^2-1}}2$.\\
Thus, we have, for $n\ge\tau$,
\bqn
D_n=\la^\tau(a\la_1^{n-\tau}+b\la_2^{n-\tau}),
\eqn
where $a,b$ can be determined by $D_\tau$ and $D_{\tau+1}$, i.e.
\bqn
\la^\tau&=&\la^\tau(a+b)\\
\la^{\tau+1}-\frac{\la^{\tau-1}}4&=&\la^\tau(a\la_1+b\la_2).
\eqn
Substitute $a,b$ into the equation $D_n=0$ and use the facts that $\la_1+\la_2=\la$ and $\la_1\la_2=\frac14$, we have
$$\frac{\la^{\tau-1}(\la_1^{n-\tau+2}-\la_2^{n-\tau+2})}{\la_1-\la_2}=0.$$
Therefore, if $\la\ne0$, we must have
$$\Big(\frac{\la_2}{\la_1}\Big)^{n-\tau+2}=1,\quad \la_1\ne\la_2,$$
from which we obtain $\frac{\la_2}{\la_1}=\cos\frac{2k\pi}{n-\tau+2}+i\sin\frac{2k\pi}{n-\tau+2}$, $k=1,2,\cdots,n-\tau+1$ ($k=0$ corresponds to the case $\la_1=\la_2$ and thus is rejected).\\
Hence, among the $n$ eigenvalues of $\textbf{C}_{n,\tau}$, $\tau-1$ of them are 0 and the rest $n-\tau+1$ are
$$\la=\la_1+\la_2=\cos\frac{k\pi}{n-\tau+2},\quad k=1,2,\cdots,n-\tau+1.$$
Define a uniform random variable $K$ taking values in $\{1,2,\cdots,n-\tau+1\}$. Then we have
\bqn
\rP(\la\le t)&=&\frac{\tau-1}nI_{[0,\infty)}(t)+\frac{n-\tau+1}n\rP\Big(\frac{K}{n-\tau+2}\ge\frac{\cos^{-1}(t)}{\pi}\Big)\\
&\to&1-\frac{\cos^{-1}(t)}{\pi}=:H(t), \quad t\in(-1,1)
\eqn
since $\frac{K}{n-\tau+2}\overset{D}{\rightarrow} \textrm{Uniform}(0,1)$ as $n\to\infty$.\\
Taking the derivative, we have  $$H'(t)=\frac{1}{\pi\sqrt{1-t^2}}, \quad t\in(-1,1).$$
The proof of the lemma is complete.
\section{Derivation of the Stieltjes Transform}
To derive the Stieltjes transform, we mainly use Theorem 4.1 in Bai and Silverstein (2010).
\begin{thm}(Theorem  4.1 in Bai and Silverstein (2010)) \label{thm2} Suppose that the entries of $X_n$ $(p\times n)$ are independent complex random variables satisfying
\bqa \label{trun2}
 \frac{1}{\eta^2np}\sum_{jk}\rE(|x_{ij}^{(n)}|^2I(|x_{ij}^{(n)}|\geq \eta\sqrt{n}))\to0.
\eqa and that $T_n$ is a sequence of Hermitian matrices independent of $X_n$ and that
the empirical spectral distribution (ESD) of $T_n$ tends to a non-random limit $H$ in some sense
(in probability or a.s.). If $p/n\to y\in(0,\infty)$, then the ESD of the product $S_nT_n$ tends to a
nonrandom limit $F$ in probability or almost surely (accordingly), where $S_n=\frac1nX_nX_n^*$.
\end{thm}
\begin{remark} Note that the eigenvalues of the product matrix $S_nT_n$ are all real although it is not
symmetric, because the whole set of eigenvalues is the same as that of the symmetric matrix $S_n^{1/2}T_nS_n^{1/2}$.
\end{remark}
\begin{remark} Note that condition (\ref{trun2}) can be implied by condition (\ref{trun1}), so Theorem \ref{thm2} is applicable to our case. In addition, the $(2+\delta)$-th moment assumption can be weakened to the 2nd moment condition.
\end{remark}
Also, according to (4.4.4) in Bai and Silverstein (2010),
\bqa\label{eq2}
\frac1zm_F\Big(\frac1z\Big)=\frac 1y-1+\frac1{2\pi iyz}\oint_{|\zeta|=\rho}\log\Big(1-z\zeta^{-1}+zy\zeta^{-1}+\zeta^{-2}zym_H\Big(\frac1\zeta\Big)\Big)d\zeta.
\eqa
Replacing $z^{-1}$ by $z$, we have
\bqa\label{eq3}
zm_F(z)=\frac 1y-1+\frac z{2\pi iy}\oint_{|\zeta|=\rho}\log\Big(z-\zeta^{-1}+y\zeta^{-1}+\zeta^{-2}ym_H\Big(\frac1\zeta\Big)\Big)d\zeta,
\eqa
where we have used the fact that the integral for $\log z$ with respect to $\zeta$ on the contour $|\zeta|=\rho$ is 0.\\
Next, set $\psi(u)=-\frac1u+y\int\frac t{1+tu}dH(t)$. Then (\ref{eq3}) becomes
\bqa
zm_F(z)&=&\frac 1y-1+\frac z{2\pi iy}\oint_{|\zeta|=\rho}\log(z-\psi(-\zeta))d\zeta\non
&=&\frac 1y-1-\frac z{2\pi iy}\oint_{|\zeta|=\rho}\frac{\zeta}{z-\psi(-\zeta)}d\psi(-\zeta)\non
&=&\frac 1y-1-\frac z{2\pi iy}\oint_{\mathcal{C}}\frac{\psi^{-1}(s)}{z-s}ds.
\eqa
When $\zeta$ is in the contour $|\zeta|=\rho$ with $\rho\in(0,1/\tau_0)$, where $\tau_0$ is a truncation point of eigenvalues of $T_n$ as defined in Section 4.3.1 in Bai and Silverstein (2010) and here we can take $\tau_0=1+\ep$ for some $\ep>0$ by Lemma \ref{key}, we have $\psi^{-1}(s)=-\zeta$ being bounded. Therefore, we have $s=z$ as the only pole. \\
Moreover, as the contour $\mathcal{C}$ is the image of the contour $|\zeta|=\rho$ under the map $\zeta\mapsto\psi(-\zeta)$ and note that $\zeta$ lies on a small circle enclosing the origin. Hence, $\mathcal{C}$ encloses the whole complex plane except a small region containing the origin.
Also, by Silverstein and Bai (1995), for each
$z\in \mathds{C}^+\equiv \{z\in\mathds{C}: \Im(z)>0\}$, there exists a unique solution $\xi\in \mathds{C}^+$ such that $z=\psi(-\xi)$.
By taking $\tau_0$ large enough, we have $s=z$ in the contour $\mathcal{C}$. Therefore,
\bqn
zm_F(z)=\frac 1y-1+\frac zy\psi^{-1}(z),
\eqn
or equivalently,
\bqa
z=\psi(ym_F+\frac {y-1}z).
\eqa
Note that
\bqn
\psi(u)&=&-\frac1u+y\int\frac t{1+tu}dH(t)\\
&=&\frac{y-1}u-\frac y{\pi u}\int_{-1}^{1}\frac 1{(1+tu)\sqrt{1-t^2}}dt\\
&=&\frac{y-1}u-\frac y{2\pi u}\int_{0}^{2\pi}\frac1{1+u\cos\theta}d\theta\\
&=&\frac{y-1}u-\frac y{2\pi ui}\oint_{|s|=1}\frac2{us^2+2s+u}ds.
\eqn
The integrand has two poles at $s_1=\frac{-1+\sqrt{1-u^2}}{u}$ and $s_2=\frac{-1-\sqrt{1-u^2}}{u}$.
As $s_1s_2=1$, we must have one of them is inside the contour and the other is outside. Therefore, we have
\bqn
\psi(u)&=&\frac{y-1}u-\frac y{2\pi ui}\oint_{|s|=1}\frac2{us^2+2s+u}ds\\
&=&\frac{y-1}u\pm\frac {2y}{u^2(s_1-s_2)}\\
&=&\frac{y-1}u\pm\frac y{u\sqrt{1-u^2}},
\eqn
where the choice of $+$ or $-$ sign is determined by which of $s_{1,2}$ is inside the contour. Substitute the above expression into (8), and we have
\bqn
m_F^2z^2[1-(ym_F+\frac{y-1}z)^2]=1.
\eqn
Note that in our question, $c=\lim_{N\to\infty} N/T=1/y$. Therefore, the Stieltjes transform $\widetilde{m}$ of the LSD of $\widetilde{\bbM}_N=\frac1{N}\textbf{X}_T^*\textbf{X}_T\textbf{C}_{T,\tau}$ satisfies
\bqn
\widetilde{m}^2z^2[1-(\frac {\widetilde{m}}c+\frac{1/c-1}z)^2]=1.
\eqn
Next, the Stieltjes transform $\underline{m}$ of the LSD of $\underline{\bbM}_N=\frac1{N}\textbf{X}_T\textbf{C}_{T,\tau}\textbf{X}_T^*$ satisfies
$\underline{m}=\frac {\widetilde{m}}c+\frac{1/c-1}z$ and therefore,
\bqn
(c\underline{m}-\frac{1-c}z)^2z^2(1-\underline{m}^2)=1.
\eqn
Finally, the Stieltjes transform $m$ of LSD of $\bbM_N=\frac1{T}\textbf{X}_T\textbf{C}_{T,\tau}\textbf{X}_T^*$ satisfies
\bqn
m(z)&=&\lim_{N\to\infty}\frac1N\sum_{i=1}^N\frac1{\lambda_i-z}=\lim_{N\to\infty}\frac1N\sum_{i=1}^N\frac1{\frac N{T}\underline{\lambda}_i-z}
=\lim_{N\to\infty}\frac1{cN}\sum_{i=1}^N\frac1{\underline{\lambda}_i-\frac1cz}=\frac1c\underline m(\frac zc).
\eqn
Substituting back to the above equation, we have
\bqn
(czm(z)+c-1)^2(1-c^2m^2(z))=1,
\eqn
which is the same as (\ref{eq1}).
\bibliographystyle{plain}

 KLASMOE and School of Mathematics and Statistics, Northeast Normal University, P.R.C, Changchun 130024\\
 Email: baizd@nenu.edu.cn\\

 Department of Statistics and Applied Probability, National University of Singapore, Singapore 117546\\
 Email: wangchen2351@gmail.com
\end{document}